\documentclass{article}
\input epsf
\usepackage{graphicx}
\usepackage{amsmath,amsxtra,amssymb,latexsym, amscd}

\newtheorem{thm}{Theorem}[section]

\newtheorem{cor}[thm]{Corollary}

\begin{document}

\title{A sufficient condition for  a set of calibrated surfaces to be area-minimizing under
diffeomorphisms}

\author{\bf Doan The Hieu \thanks{This work was completed at SNU, with the support of Korean Foundation for Advanced Studies.
The author was also supported in part by the National Basic Research Program in Natural Sciences, Vietnam.}
\\
 Department of Mathematics\\
 College of Education, Hue University\\
 32 Le Loi, Hue, Viet Nam\\
\\
 dthehieu@yahoo.com; deltic@dng.vnn.vn}

\maketitle
\begin{abstract}
We extend Choe's idea in \cite{choe} to nonpolyhedral calibrated surfaces and give some examples of polyhedral sets over right prisms and  nonpolyhedral calibrated surfaces.
\end{abstract}

\maketitle

\section{Introduction}

 In \cite{choe}, Choe proved ``{\sl Every stationary polyhedral set is area-minimizing under diffeomophisms leaving the boundary fixed}".  In his proof, a system of differential forms and orientations (of faces) was chosen at each singular edge. In fact, the differential forms are calibrations that calibrate the faces at each singular edge and have the vanishing sum. We observe that, the suitable orientations of faces at each singular edge determine the same orientation on it whenever it lies on the boundary of faces.

By the above observation, we extend Choe's idea by proving a sufficient condition for   certain sets
of calibrated surfaces (including polyhedral sets) to be
area-minimizing under diffeomophisms leaving the boundary fixed.  This sufficient condition, when applies to polyhedral sets, is also necessary.

We give some more examples of polyhedral sets over right prisms and first examples of  nonpolyhedral calibrated surfaces ($2$-dimensional  ones with singular sets of dimension $1$ in $\Bbb R^4).$

The author would like to thank the gracious host scholar,
Professor Hong-Jong Kim,  for his thoughtfulness and assistance.
Many thanks to Professor Jaigyoung Choe for his help, suggestions
and the invitation talking this result at his seminar. This paper
was written while the author was visiting RIM-GARC, Department of
Mathematics of the Seoul National University, Korea. We would like
to thank that institution  for its hospitality and generous
support.

\medskip
\section{The theorem}
\bigskip

We refer the readers to \cite{choe} for the definition of polyhedral sets.

Let $\{C_i\}_{i\in I}$ be a set of calibrated surfaces of
dimension $m$ in $\Bbb R^n  (m<n)$ and $\{w_i\}_{i\in I}$ be the
set of correspondent calibrations. That means for each $i\in I,
w_i$ calibrates $C_i$ with a suitable orientation.
Note that if $\omega_i$ calibrates $C_i$, then $-\omega_i$ calibrates $C_i$ with opposite orientation. Depending on a chosen orientation on $C_i$ we have the corespondent calibration to be $\omega_i$ or $-\omega_i.$

Let  $\Sigma\subset\Bbb R^n$ be a set satisfies the following conditions:

(i)\ $\Sigma\subset \cup_{i\in I}C_i,$

(ii)\ the set $E =\Sigma\cap(C_i\cap C_j)$ is of dimension $m-1$ for every $i,j\in I, i\not =j.$

We call each $F_i=\Sigma\cap C_i$ a face, each $E$ a singular
edge, the union of all singular edges $ E$ the singular set $S$,
the closure of $\partial F_i\sim S$ the boundary edge of $\Sigma$
in $F_i,$ the union $\cup_{i\in I}( \partial F_i\sim S)$ the boundary $\partial \Sigma$ of $\Sigma.$

$\Sigma$   is said to be
area-minimizing under diffeomorphisms leaving the boundary fixed if
 $$Vol(\Sigma)\le
Vol(\varphi(\Sigma)),$$ for any diffeomorphism $\varphi$ of $\Bbb
R^n$ leaving the boundary of $\Sigma$ fixed.

Suppose $\{E_j\}_{j\in J}$ is the set of all singular edges and $\{F_i\}_{i\in I}$ is the set of all faces of $\Sigma.$  Denote
 $$I_{E_j}= \{i : F_i\supset E_j \}\subset I,$$
 $$J_{F_i}=\{j :  E_j\subset F_i \}\subset J.$$

\begin{thm}\label{t:main}
Let $\Sigma$ be a set defined as above. Suppose that every singular edge $E_j$ lies on the boundary $\partial F_i, \forall i\in I_{E_j}$ and for each $E_j$ we can choose suitable orientations on $F_i, \forall i\in I_{E_j},$ such that:

(i) the  orientations on $F_i, \forall i\in I_{E_j}$ determine the same orientation on $E_j,$

(ii) the corespondent calibrations have vanishing sum.

 \noindent Then $\Sigma$ is
area-minimizing under diffeomorphisms leaving $\partial\Sigma$
fixed.
\end{thm}

\emph{Proof.}
The reasonings of the proof are very similar as that of the main theorem in \cite{choe} with some little changes.

Let $\varphi$ be a diffeomorphism leaving $\partial
\Sigma$ fixed and $\varphi_t$ be the homotopy from the identity
to $\varphi.$ Suppose $G_j$ is the $m$-dimensional smooth surface
swept out by $\varphi_t(E_j)$ and $D_i$ is $(m+1)$-dimensional
surface swept out by $\varphi_t(F_i).$ We have
$$\partial D_i=F_i\cup\varphi(F_i)\cup_{j\in J_{F_i}}G_j,$$
and hence
$$\int_{\partial
D_i}w_i=\int_{F_i}w_i+\int_{\varphi(F_i)}w_i+\sum_{j\in
J_{F_i}}\int_{G_j}w_i.$$

  Since $w_i$ is a
calibration that calibrates $F_i$,  we get the following inequality:
$$Vol(F_i)\le Vol(\varphi(F_i))-\sum_{j\in J_{F_i}}\int_{G_j}w_i,$$
and finally
$$Vol(\Sigma)\le Vol(\varphi(\Sigma))-\sum_{j\in J}\sum_{i\in
I_{E_j}}\int_{G_j}w_i.$$

By virtue of the assumtions of the theorem, we can assume the  orientations  on  $F_i, \forall i\in I_{E_j},$ determine the same orientation on $G_j$ and since $\sum_{i\in
I_{E_j}}w_i=0,$ the last term equals  zero.
The theorem is proved.

\begin{cor}\label{t:main2}
Let $\Sigma$ be a polyhedral set.  Then $\Sigma$ is
area-minimizing under diffeomorphisms leaving $\partial\Sigma$
fixed if and only if  $\Sigma$   satisfies the assumptions in the Theorem    ~\ref{t:main}.
\end{cor}
\emph{Proof.}  The sufficiency  follows from  the above theorem and the necessity  follows from the proof of the main theorem in \cite{choe}.

\section{Examples}

1. At Ken Brakke's homepage,
    http://www.susqu.edu/facstaff/b/brakke/, we can see eight  nice polyhedral cones, that are made of flat sheets meeting along
   triple lines with an equal angle $120^0.$  All of them are  area-minimizing under diffeomophisms leaving the boundary fixed by virtue of Theorem ~\ref{t:main}.
 Figures 1 provide more three polyhedral sets over  right prisms.
\begin{figure}[htb]
\includegraphics{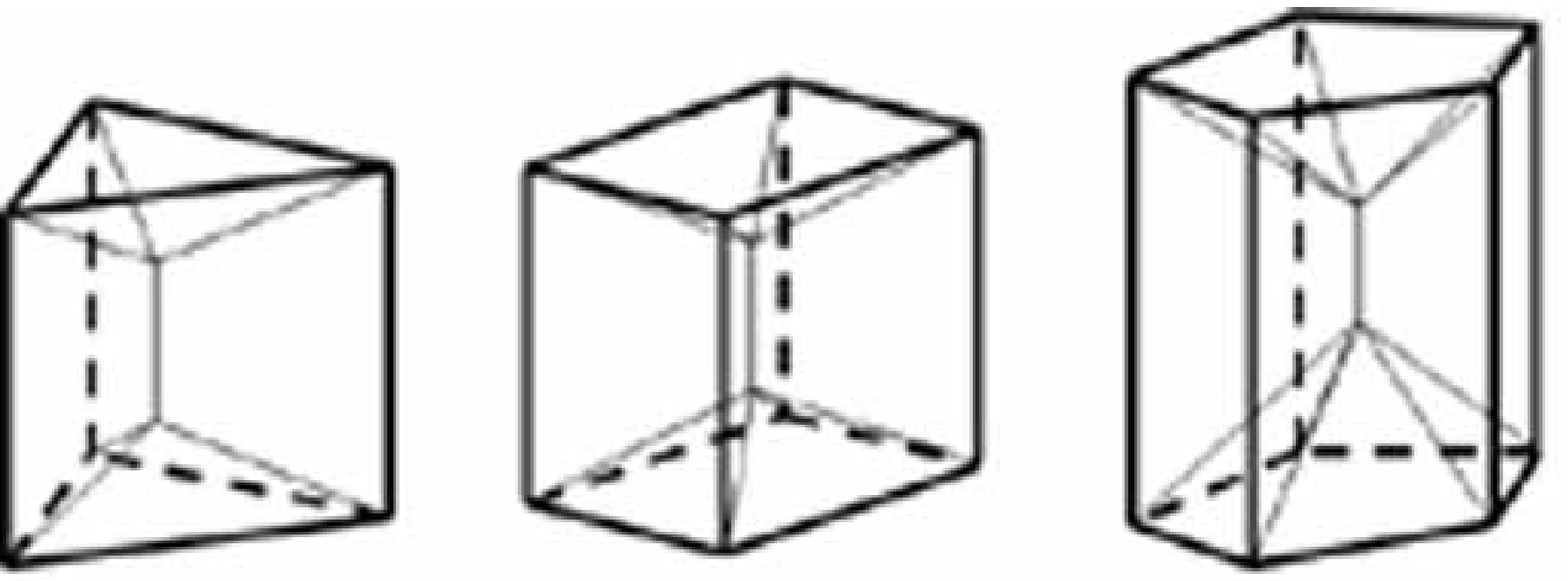}
\caption{Polyhedral sets over right prisms.}
\end{figure}

2. Below are examples of nonpolyhedral calibrated surfaces that is area-minimizing under diffeomophisms leaving the boundary fixed.

Let $\Bbb C^2\equiv \Bbb R^4$ be complex plane with the
standard complex structure $J_1,$ \ \ $J_1e_1=e_3; J_1e_2=e_4.$

Let $R_2, R_3,\ldots, R_n$ be the rotations of angles $\alpha, 2\alpha, \ldots, (n-1)\alpha$
about the plane $\{x_3=x_4= 0\},$ respectively, where $\alpha$ satisfies the condition $n\alpha=2\pi, n\in \Bbb N.$ And let $J_2, J_3, \ldots, J_n$ be $(n-1)$ complex structures  on $\Bbb R^4$
 induced by $R_2, R_3,\ldots, R_n;$
      $$J_i(e_1)=R_i(e_3),\ \
J_i(e_2)=R_i(e_4);\ \  i=2,3,\ldots, n.$$

Denote $w_1, w_2,\ldots, w_n$ the K\"{a}hler forms correspondent
to $J_1, J_2,\ldots, J_n .$
We can easily to see that:
$$\sum^n_{i=1} w_i=0.$$

 Consider the complex curves:
$$\begin{aligned}
C&=\{(z,w)\in \Bbb C^2 : z=w^2\}\\
&=\{(x_1, x_2, x_3, x_4) : \ \ x_2=x_1^2-x_3^2;\ \ x_4=2x_1x_3\}.\\
\end{aligned}$$  

Let $D$ be the intersection of $C$ and  $\{ \sum_{i=1}^4x_i^2=1;  x_1\ge 0; x_3\ge 0 \}.$ Note that $D$ contains two  planar curves
$\{x_2=x_1^2; \ x_3=x_4=0; x_1^2+x_2^2\le 1; x_1\ge 0\},$
 and
$\{x_2=-x_3^2; \ x_1=x_4=0; x_2^2+x_3^2\le 1; x_3\ge 0\}.$ By using  the rotations of angles $k\alpha,\ k=1,2.\ldots, n-1$
about the plane $\{x_3=x_4= 0\}  $  we get the images $D_{i}\  (i=2,3,\ldots, n)$ of $D.$
  Obviously, $w_1$ calibrates $D$ and $w_i$ calibrates $D_i,\ i=2,3, \ldots n.$

The set $\Sigma= D\bigcup D_i$  contains one singular edge and is area-minimizing under diffeomophisms leaving the boundary fixed by virtue of Theorem ~\ref{t:main}.

Similarly, by using  the rotations of angles $k\beta,\ k=1,2.\ldots, m-1; m\beta=2\pi$
about the plane $\{x_1=x_4= 0\} $, we get the images $D'_{j}\ (j=2,3,\ldots, m)$ of $D$
  and the images $\Sigma_{j}\  (j=2,3,\ldots, m)$  of $\Sigma.$

 The set $D\bigcup D_i\bigcup D'_j$ contains two singular edges. The set $ \Sigma'= \Sigma\bigcup\Sigma_j$ contains many singular edges. By the same reasoning as above, they are also area-minimizing under diffeomophisms leaving the boundary fixed.

\end{document}